\documentclass[11pt]{amsart}
\usepackage{amsmath,amssymb,amsfonts,url,mathptmx}
\numberwithin{equation}{section}

\newtheorem{thm}{Theorem}[section]
\newtheorem{lem}{Lemma}[section]
\newtheorem{rem}{Remark}[section]
\newtheorem{prop}{Proposition}[section]

\begin{document}
\title[Singular Liouville Equation]{Vanishing Estimates for Liouville equation with quantized singularities}

\author{Juncheng Wei}
\address{Department of Mathematics \\ University of British Columbia\\ Vancouver, BC V6T1Z2, Canada} \email{jcwei@math.ubc.ca }

\author{Lei Zhang} \footnote{The research of J. Wei is partially supported by NSERC of Canada. Lei Zhang is partially supported by a Simons Foundation Collaboration Grant}
\address{Department of Mathematics\\
        University of Florida\\
        1400 Stadium Rd\\
        Gainesville FL 32611}
\email{leizhang@ufl.edu}

\date{\today}

%%%%%%%%%%%%%%%%%%%%%%%%%%%%%%%%%%%%%%%%%%%%%
\begin{abstract} In this article we continue with the research initiated in our previous work on singular Liouville equations with quantized singularity. The main goal of this article is to prove
that as long as the bubbling solutions violate the spherical Harnack inequality near a singular source, the first derivatives of coefficient functions must tend to zero.
\end{abstract}

%%%%%%%%%%%%%%%%%%%%%%%%%%%%%%%%%%%%%%%%%%%%%
\maketitle

\section{Introduction}
In this article we study bubbling solutions of
\begin{equation}\label{main-eq}
\Delta u+ \mathrm{H}(x) e^{u}=4\pi \alpha \delta_0\quad \mbox{ in }\quad \Omega\subset \mathbb R^2
\end{equation}
where $\Omega$ is an open, bounded subset of $\mathbb R^2$ that contains the origin, $\alpha>-1$ is a constant and $\delta_0$ is the Dirac mass at $0$, $\mathrm{H}$ is a positive and smooth function. Since
$$\Delta (\frac 1{2\pi}\log |x|)=\delta_0, $$
 one can use the logarithmic function $2\pi \alpha\log |x|$ to remove the singular source from the equation: let $u_1(x)=u(x)-2\alpha \log |x|$, then $u_1$ satisfies

\begin{equation}\label{main-2}
\Delta u_1+|x|^{2\alpha}\mathrm{H}(x)e^{u_1}=0,\quad \mbox{in}\quad \Omega.
\end{equation}

 If a sequence of solutions $\{u^k\}_{k=1}^{\infty}$ of (\ref{main-2}) satisfies
 $$\lim_{k\to \infty}u^k(x_k)=\infty,\quad \mbox{ for some $\bar x\in B_{\tau}$ and $x_k\to \bar x$.} $$
  we say $u^k$ is a sequence of bubbling solutions or blowup solutions, $\bar x$ is called a blowup point.  For many reasons in applications it is most interesting to consider $\alpha=N\in \mathbb N$ (the set of natural numbers) and when $0$ is the only blowup point of $u^k$. Our set-up of bubbling solutions is as follows:
Let $\mathfrak{u}_k$ be a sequence of solutions of
\begin{equation}\label{t-u-k}
\Delta \mathfrak{u}_k(x)+|x|^{2N}\mathrm{H}_k(x)e^{\mathfrak{u}_k}=0, \quad \mbox{in}\quad B_{\tau}
\end{equation}
for some $\tau>0$ independent of $k$. $B_{\tau}$ is the ball centered at the origin with radius $\tau$.  In addition we postulate the usual assumptions on $\mathfrak{u}_k$ and $\mathrm{H}_k$:
For a positive constant $C$ independent of $k$, the following holds:
\begin{equation}\label{assumption-1}
\left\{\begin{array}{ll}
\|\mathrm{H}_k\|_{C^3(\bar B_{\tau})}\le C, \quad \frac 1C\le \mathrm{H}_k(x)\le C, \quad x\in \bar B_{\tau}, \\ \\
\int_{B_{\tau}} \mathrm{H}_k e^{\mathfrak{u}_k}\le C,\\  \\
|\mathfrak{u}_k(x)- \mathfrak{u}_k(y)|\le C, \quad \forall x,y\in \partial B_{\tau},
\end{array}
\right.
\end{equation}
and since we study the asymptotic behavior of blowup solutions around the singular source, we assume that there is no blowup point except at the origin:
\begin{equation}\label{assump-2}
\max_{K\subset\subset B_{\tau}\setminus \{0\}} \mathfrak{u}_k\le C(K).
\end{equation}
Also,
we use the value of $\mathfrak{u}_k$ on $\partial B_{\tau}$ to define a harmonic function $\phi_k(x)$:
\begin{equation}\label{phi-k}
\left\{\begin{array}{ll}
\Delta \phi_k(x)=0,\quad \mbox{in}\quad B_{\tau},\\ \\
\phi_k(x)=\mathfrak{u}_k(x)-\frac{1}{2\pi \tau}\int_{\partial B_{\tau}}\mathfrak{u}_kdS,\quad x\in \partial B_{\tau}.
\end{array}
\right.
\end{equation}
Clearly the mean value property of harmonic functions implies $\phi_k(0)=0$ and the finite oscillation of $\mathfrak{u}_k$ on $\partial B_{\tau}$ means that all derivatives of $\phi_k$ are uniformly bounded in $B_{\tau/2}$.
In this article we consider the case that:
\begin{equation}\label{no-sp-h}
\max_{x\in B_{\tau}} \mathfrak{u}_k(x)+2(1+N)\log |x|\to \infty,
\end{equation}
which is equivalent to saying that the spherical Harnack inequality does not hold for $\mathfrak{u}_k$. It is also mentioned in literature ( see \cite{kuo-lin-jdg, wei-zhang-adv} ) that $0$ is called an non-simple blowup point. The main result of this article is

\begin{thm}\label{main-thm} Let $\{\mathfrak{u}_k\}$ be a sequence of solutions of (\ref{t-u-k}) such that (\ref{assumption-1}),(\ref{assump-2}) and (\ref{no-sp-h}) hold. Then
$$\nabla (\log\mathrm{H}_k+\phi_k)(0)=o(1),\quad \mbox{as } k\to \infty, $$ where $\phi_k$ is defined in (\ref{phi-k}).
\end{thm}

When bubbling solutions satisfy (\ref{no-sp-h}), they are called non-simple blowup solutions ( see \cite{kuo-lin-jdg}). Theorem \ref{main-thm} is a complement of Theorem 1.1 of \cite{wei-zhang-adv}, which asserts that under certain conditions (see Theorem A below) $\nabla(\log \mathrm{H}_k+\phi_i^k)(0)$ tends to zero. Theorem \ref{main-thm} removes the restrictions in \cite{wei-zhang-adv}. In other words, the combination of Theorem A and Theorem \ref{main-thm} proves that $\nabla (\log \mathrm{H}_k+\phi_k)(0)\to 0$ as long as the non-simple blowup situation occurs. Besides the advancement of analytical understanding, this conclusion is particularly important in application.
Theorem \ref{main-thm} can be applied to situations beyond single equations. For certain systems of equations such as Toda systems, the bubble accumulations can be described by a sequence of bubbling solutions with quantized singular source. Theorem \ref{main-thm} is very useful to rule out complicated bubbling accumulation pictures in Toda systems.

It remains an open question whether or not $\nabla (\log\mathrm{H}_k+\phi_k)(0)$ tends to $0$ if the bubbling solutions satisfy the spherical Harnack inequality around the origin. We tend to believe one can construct a sequence of bubbling solutions that satisfy spherical Harnack inequality with nonzero first derivatives of the coefficients functions. In particular the works of Del Pino-Esposito-Musso \cite{del-pino-1, del-pino-2} on two dimensional Euler flows seem to suggest that for simple blowup solutions with quantized singular sources, the first derivatives of the coefficient functions may not tend to zero at singular sources. Instead the $ (N+1)$-th derivatives  should vanish. Our result, together with \cite{del-pino-1, del-pino-2},  demonstrates a  striking contrast  between simple and non-simple bubblings.

The non-simple bubbling situation and vanishing theorems have profound impact to problems in geometry and physics. For example
for the following mean field equation defined on a Riemann surface $(M,g)$:
\begin{equation}\label{mean-sin}
\Delta_gu+\rho(\frac{h(x)e^{u(x)}}{\int_Mhe^{u}}-\frac 1{Vol_g(M)})=4\pi\sum_j \alpha_j (\delta_{p_j}-\frac 1{Vol_g(M)}),
\end{equation}
the solution $u$ represents a conformal metric with prescribed conic singularities (see \cite{erem-3,tr-1,tr-2}) . In particular if the singular source is quantized, the Liouville equation has close ties with Algebraic geometry, integrable system, number theory and complex Monge-Ampere equations (see \cite{chen-lin-last-cpam}). In Physics the understanding of non-simple blowup phenomenon would be extremely useful for the study of mean field limits of point vortices in the Euler flow \cite{caglioti-1,caglioti-2} and models in the Chern-Simons-Higgs theory \cite{jackiw} and in the electroweak theory \cite{ambjorn}, etc. It is also remarkable that non-simple bubbling solutions also occur in systems.
In \cite{gu-zhang}, the non-simple blowup solutions are studied for singular Liouville systems. Finally we remark that when the blowup point is a location of a singular source, whether or not this point has to be a critical point of coefficient functions has intrigued people for years. Our previous result \cite{wei-zhang-adv} is the first result for singular Liouville equation, the second author proved a surprising vanishing theorem for singular Toda systems in \cite{zhang-imrn}.

The organization of this paper is as follows. In section two we review a few fundamental tools for the proof of the main theorem and invoke several key estimates established in our previous work \cite{wei-zhang-adv}. Then in section three we use a sequence of global solutions to approximate our blowup solutions. The point-wise estimates proved in this section are more precise than what is established in \cite{wei-zhang-adv} and are important for our argument. In section four we prove a crucial estimate on the difference between blowup solution and the global solutions as the first term in the approximation. As a consequence of section four, we move to section five to complete the proof of the main theorem. The proof in section four is similar to the proof of uniqueness theorems for bubbling solutions in \cite{wu-zhang-1}, \cite{yang-wen-1}, \cite{lin-yan-uniq}, etc.

\medskip

{\bf Notation:} We will use $B(x_0,r)$ to denote a ball centered at $x_0$ with radius $r$. If $x_0$ is the origin we use $B_r$. $C$ represents a positive constant that may change from place to place.

\section{Preliminary discussions}

For simple notation we set
\begin{equation}\label{uk-d}
u_k(x)=\mathfrak{u}_k(x)-\phi_k(x), \quad \mbox{and}
\end{equation}
\begin{equation}\label{hk-d}
h_k(x)=\mathrm{H}_k(x)e^{\phi_k(x)}
\end{equation}
to write the equation of $u_k$ as
\begin{equation}\label{eq-uk}
\Delta u_k(x)+|x|^{2N}h_k(x)e^{u_k}=0,\quad \mbox{ in }\quad B_{\tau}
\end{equation}
Without loss of generality we assume
\begin{equation}\label{rea-h}
\lim_{k\to \infty} h_k(0)=1.
\end{equation}

Obviously (\ref{no-sp-h}) is equivalent to
\begin{equation}\label{no-sp-h-u}
\max_{x\in B_{\tau}} u_k(x)+2(1+N)\log |x|\to \infty.
\end{equation}

It is well known \cite{kuo-lin-jdg, bart3} that $ u_k$ exhibits a non-simple blowup profile.  It is established in \cite{kuo-lin-jdg,bart3} that there are $N+1$ local maximum points of $ u_k$: $p_0^k$,....,$p_{N}^k$ and they are evenly distributed on $\mathbb S^1$ after scaling according to their magnitude: Suppose along a subsequence
$$\lim_{k\to \infty}p_0^k/|p_0^k|=e^{i\theta_0}, $$
then
$$\lim_{k\to \infty} \frac{p_l^k}{|p_0^k|}=e^{i(\theta_0+\frac{2\pi l}{N+1})}, \quad l=1,...,N. $$
For many reasons it is convenient to denote $|p_0^k|$ as $\delta_k$ and define $\mu_k$ as follows:
\begin{equation}\label{muk-dk}
\delta_k=|p_0^k|\quad \mbox{and }\quad \mu_k= u_k(p_0^k)+2(1+N)\log \delta_k.
\end{equation}

Since $p_l^k$'s are evenly distributed
around $\partial B_{\delta_k}$, standard results for Liouville equations around a regular blowup point can be applied to have $ u_k(p_l^k)= u_k(p_0^k)+o(1)$. Also, (\ref{no-sp-h}) gives $\mu_k\to \infty$. The interested readers may look into \cite{kuo-lin-jdg,bart3} for more detailed information.

In our previous work \cite{wei-zhang-adv} we prove the following vanishing type estimates for the first derivatives of the coefficient function $\log h_k$:

\medskip

\emph{Theorem A: Let $u_k$, $\phi_k$, $h_k$, $\delta_k$, $\mu_k$ be defined by (\ref{eq-uk}), (\ref{phi-k}), (\ref{hk-d}), (\ref{muk-dk}) respectively. Then
\begin{equation}\label{first-order-v}
|\nabla \log h_k(0)|=O(\delta_k)+O(\delta_k^{-1}e^{-\mu_k}\mu_k).
\end{equation}
}

Here we observe that if $\mu_ke^{-\mu_k}=o(\delta_k)$, we already have $\nabla h_k(0)=o(1)$, which is equivalent to $\nabla (\mathrm{H}_ke^{\phi_k})(0)=o(1)$. Thus throughout the paper we assume
\begin{equation}\label{small-delta}
\delta_k\le C\mu_ke^{-\mu_k}.
\end{equation}
Finally we shall use $E$ to denote a frequently appearing error term of the size $O(\delta_k^2)+O(\mu_ke^{-\mu_k})$. Because of (\ref{small-delta}), $$E=O(\mu_ke^{-\mu_k}).$$

\section{Approximating bubbling solutions by global solutions}

First we recall that $|p_0^k|=\delta_k$, so we write $p_0^k$ as $p_0^k=\delta_ke^{i\theta_k}$ and define $v_k$ as
\begin{equation}\label{v-k-d}
v_k(y)=u_k(\delta_k ye^{i\theta_k})+2(N+1)\log \delta_k,\quad |y|<\tau \delta_k^{-1}.
\end{equation}
If we write out each component, (\ref{v-k-d}) is
$$
v_k(y_1,y_2)=u_k(\delta_k(y_1\cos\theta_k-y_2\sin\theta_k),\delta_k(y_1\sin\theta_k+y_2\cos\theta_k))+2(1+N)\log \delta_k. $$
Then it is standard to verify that $v_k$ solves

\begin{equation}\label{e-f-vk}
\Delta v_k(y)+|y|^{2N}\mathfrak{h}_k(\delta_k y)e^{v_k(y)}=0,\quad |y|<\tau/\delta_k,
\end{equation}
where
\begin{equation}\label{frak-h}
\mathfrak{h}_k(x)=h_k(xe^{i\theta_k}),\quad |x|<\tau.
\end{equation}
Thus the image of $p_0^k$ after scaling is $Q_1^k=e_1=(1,0)$.
Let $Q_1^k$, $Q_2^k$,...,$Q_{N}^k$ be the images of $p_i^k$ $(i=1,...,N)$ after the scaling:
$$Q_l^k=\frac{p_l^k}{\delta_k}e^{-i\theta_k},\quad l=1,...,N. $$
 It is established by Kuo-Lin in \cite{kuo-lin-jdg} and independently by Bartolucci-Tarantello in \cite{bart3} that
\begin{equation}\label{limit-q}
\lim_{k\to \infty} Q_l^k=\lim_{k\to \infty}p_l^k/\delta_k=e^{\frac{2l\pi i}{N+1}},\quad l=0,....,N.
\end{equation}
Then in our previous work \cite{wei-zhang-adv} we obtained ( see (3.13) in \cite{wei-zhang-adv})
\begin{equation}\label{Qm-close}
Q_l^k-e^{\frac{2\pi l i}{N+1}}=E.
\end{equation}
Choosing $3\epsilon>0$ small and independent of $k$, we can make disks centered at $Q_l^k$ with radius $3\epsilon$ (denoted as $B(Q_l^k,3\epsilon ) $) mutually disjoint. Let
\begin{equation}\label{v-muk}
\mu_k=\max_{B(Q_0^k,\epsilon)} v_k.
\end{equation}
Since $Q_l^k$ are evenly distributed around $\partial B_1$, it is easy to use standard estimates for single Liouville equations (\cite{zhangcmp,gluck,chenlin1}) to obtain
$$\max_{B(Q_l^k,\epsilon)}v_k=\mu_k+o(1),\quad l=1,...,N. $$

Let
\begin{equation}\label{def-Vk}
V_k(x)=\log \frac{e^{\mu_k}}{(1+\frac{e^{\mu_k}\mathfrak{h}_k(\delta_k e_1)}{8(1+N)^2}|y^{N+1}-e_1|^2)^2}.
\end{equation}
Clearly $V_k$ is a solution of
\begin{equation}\label{eq-for-Vk}
\Delta V_k+\mathfrak{h}_k(\delta_k e_1)|y|^{2N}e^{V_k}=0,\quad \mbox{in}\quad \mathbb R^2, \quad V_k(e_1)=\mu_k.
\end{equation}
This expression is based on the classification theorem of Prajapat-Tarantello \cite{prajapat}.

The estimate of $v_k(x)-V_k(x)$ is important for the main theorem of this article. For convenience we use
$$\beta_l=\frac{2\pi l}{N+1}, \quad \mbox{so}\,\, e_1=e^{i\beta_0}=Q_0^k,\quad
e^{i\beta_l}=Q_l^k+E,\,\,\mbox{ for }\,\, l=1,...,N. $$

\begin{prop}\label{vk-Vk} Let $l=0,...,N$ and $\delta$ be small so that $B(e^{i\beta_l},\delta)\cap B(e^{i\beta_s},\delta)=\emptyset$ for $l\neq s$.
In each $B(e^{i\beta_l},\delta)$
\begin{equation}\label{global-close}
|v_k(x)-V_k(x)|\le \left\{\begin{array}{ll}
C\mu_ke^{-\mu_k/2},\quad |x-e^{i\beta_l}|\le Ce^{-\mu_k/2}, \\
\\
C\frac{\mu_ke^{-\mu_k}}{|x-e^{i\beta_l}|}+O(\mu_k^2e^{-\mu_k}),\quad Ce^{-\mu_k/2}\le |x-e^{i\beta_l}|\le \delta.
\end{array}
\right.
\end{equation}
\end{prop}

\begin{rem} Once (\ref{global-close}) is established.  We shall use a re-scaled version of Proposition \ref{vk-Vk}. Let $\epsilon_k=e^{-\frac 12\mu_k}$, we have
\begin{equation}\label{vk-Vk-2}
|v_k(e^{i\beta_l}+\epsilon_ky)-V_k(e^{i\beta_l}+\epsilon_ky)|\le C\mu_k^2 \epsilon_k (1+|y|)^{-1},\quad 0<|y|<\delta_0 \epsilon_k^{-1}.
\end{equation}
\end{rem}

\noindent{\bf Proof of Proposition \ref{vk-Vk}}: The main idea of the proof is as follows. First from the Green's representation of $v_k$ we obtain a rather precise estimate of $v_k$ in $B_3$ away from bubbling disks. On the other hand around each $Q_m^k$ we invoke a standard pointwise estimate in \cite{chenlin1,zhangcmp,gluck} for Liouville equation around a blowup point, which provides a precise description of $v_k$ in a neighborhood of $Q_m^k$. The comparison of these two estimates gives an accurate estimate of the maximum of $v_k$ around each local maximum point.

Fixing the neighborhood of one $Q_m^k$, we first cite a result of Gluck \cite{gluck} (Appendix B of \cite{wei-zhang-adv}) to write $v_k$ as
\begin{equation}\label{tem-vk-2}
v_k(y)=\log \frac{e^{\mu_{k,m}}}{(1+e^{\mu_{k,m}}\frac{|\tilde Q_m^k|^{2N}\mathfrak{h}_k(\delta_k\tilde Q_m^k)}8|y-\tilde Q_m^k|^2)^2}+\phi_m^k(y)+O(\mu_k^2e^{-\mu_k})
\end{equation}
where $\mu_{k,m}= v_k(\tilde Q_m^k)$, $\phi_m^k$ is the harmonic function taking $0$ at $Q_m^k$ that makes $v_k-\phi_m^k=constant$ on $\partial B(Q_m^k,\delta)$.  $\tilde Q_m^k$ is where $v_k-\phi_m^k$ takes its local maximum in a neighborhood of $Q_m^k$. The difference between $\tilde Q_m^k$ and $Q_m^k$ is $O(e^{-\mu_k})$.  First we claim that
\begin{equation}\label{mu-lk}
\mu_{k,m}-\mu_k=E.
\end{equation}

Let $\Omega_k=B(0,\tau \delta_k^{-1})$ and $G_k(y,\eta)$ be the Green's function on $\Omega_k$ with respect to Dirichlet boundary condition:
\begin{equation}\label{expre-G}
G_k(y,\eta)=-\frac 1{2\pi}\log |y-\eta |+\frac{1}{2\pi}\log (\frac{|\eta |}{\tau \delta_k^{-1}}|\frac{\tau^2\delta_k^{-2}\eta}{|\eta |^2}-y|),
\end{equation}
 then for $y$ away from bubbling areas and $|y|\sim 1$, we have
\begin{align*}
v_k(y)&=v_k|_{\partial \Omega_k}+\int_{\Omega_k}G_k(y,\eta)\mathfrak{h}_k(\eta)|\eta |^{2N}e^{v_k}d\eta, \\
&=v_k|_{\partial \Omega_k}+\sum_{l=0}^NG(y,Q_l^k)\int_{B(Q_l^k,\epsilon)}|\eta |^{2N}\mathfrak{h}_k(\delta_k\eta)e^{v_k}d\eta \\ &+\sum_{l}\int_{B(Q_l^k,\epsilon)}(G_k(y,\eta)-G_k(y,Q_l^k))|\eta|^{2N}\mathfrak{h}_k(\delta_k\eta)e^{v_k}d\eta+E,\\
&=v_k|_{\partial \Omega_k}+8\pi \sum_l G_k(y,Q_l^k)+E.
\end{align*}
Later we shall use $H_k$ to denote the second term in (\ref{expre-G}). Note that we use two standard estimates. First the integration outside bubbling disk is $E$ because
$$v_k(x)\le -\mu_k-(4(N+1)-o(1))\log |x|+C, \quad 3<|x|<\tau \delta_k^{-1}. $$ Second, in the evaluation of the integral terms above we use standard bubble expansion formula (see Gluck \cite{gluck}, for example) and symmetry properties. This part is mentioned in Lemma 2.1 and Appendix B of \cite{wei-zhang-adv}. It is important to point out that the second estimate does not depend on $m$.
 In particular if we consider $y$ located at $|y-Q_m^k|=\epsilon$, the expression of $v_k$ can be written as
\begin{align}\label{tem-vk}
v_k(y)&=v_k|_{\partial \Omega_k}-4\log |y-Q_m^k|+\phi_m^k\\
&-4\sum_{l=0,l\neq m}^N\log |Q_m^k-Q_l^k|+8\pi\sum_{l=0}^N H_k(Q_m^k,Q_l^k)+E,\nonumber
\end{align}
where
\begin{equation}\label{har-around-q}
\phi_m^k=\sum_{l=0,l\neq m}^N(-4)\log \frac{|y-Q_l^k|}{|Q_m^k-Q_l^k|}+8\pi\sum_{l=0}^N(H_k(y,Q_l^k)-H_k(Q_m^k,Q_l^k))
\end{equation}
is the harmonic function that takes $0$ at $Q_m^k$ and eliminates the oscillation of $v_k$ on $\partial B(Q_m^k,\epsilon)$.
On the other hand from (\ref{tem-vk-2}) we have
\begin{equation}\label{tem-vk-8}
v_k(y)=-\mu_{k,m}-2\log \frac{|\tilde Q_m^k|^{2N}\mathfrak{h}_k(\delta_k\tilde Q_m^k)}{8}-4\log |y-Q_m^k|+\phi_m^k+O(\mu_k^2e^{-\mu_k}).
\end{equation}
Comparing (\ref{tem-vk-8}) and (\ref{tem-vk}) on $|y-Q_m^k|=\epsilon$ we have
\begin{align}\label{tem-vk-3}
&-\mu_{m,k}-2\log \frac{|\tilde Q_m^k|^{2N}\mathfrak{h}_k(\delta_k\tilde Q_m^k)}{8}\\
=&-4\sum_{l=0,l\neq m}^N\log |Q_m^k-Q_l^k|+8\pi\sum_{l=0}^N H_k(Q_m^k,Q_l^k)+v_k|_{\partial \Omega_k}+O(\mu_k^2e^{-\mu_k}). \nonumber
\end{align}
To evaluate terms in (\ref{tem-vk-3}) we observe that (see (\ref{Qm-close}))
$$|\tilde Q_m^k|^{2N}=1+E,\qquad \mathfrak{h}_k(\delta_k\tilde Q_m^k)=1+E,$$
$$Q_m^k=e^{i\beta_m}+E, \qquad \tilde Q_m^k=Q_m^k+O(e^{-\mu_k}),$$
and by  the expression of $H_k(y,\eta)$:
\begin{align*}
H_k(y,\eta)=\frac 1{2\pi}\log (\frac{|\eta |}{\tau \delta_k^{-1}}|\frac{\tau^2\delta_k^{-2}\eta }{|\eta |^2}-y|)\\
=\frac 1{2\pi}\log (\tau \delta_k^{-1})+ \frac 1{2\pi}\log |\frac{\eta }{|\eta |}-\frac{|\eta |}{\tau^2}\delta_k^2y|
\end{align*}
 we have
$$H_k(Q_m^k,Q_l^k)=\frac{1}{2\pi}\log (\tau \delta_k^{-1})+E.$$
Thus two terms in (\ref{tem-vk-3}) are
\begin{equation}\label{tri-1}
8\pi\sum_{l=0}^N H_k(Q_m^k,Q_l^k)=4(N+1)\log (\tau \delta_k^{-1})+E
\end{equation}
\begin{equation}\label{tri-2}
\sum_{l=0,l\neq m}^N\log |Q_m^k-Q_l^k|=\sum_{l=0,l\neq m}^N\log |e^{i\beta_m}-e^{i\beta_l}|+E
=\log (N+1)+E.
\end{equation}
The last equality can be verified by direct computation.
Using (\ref{tri-1}) and (\ref{tri-2}) in (\ref{tem-vk-3}) we have
\begin{align}\label{vk-bry}
v_k|_{\partial \Omega_k}=&-\mu_{m,k}-2\log \frac{\mathfrak{h}_k(\delta_ke_1)}8+4\log (1+N)-4(1+N)\log (\tau \delta_k^{-1})\\
&+O(\mu_k^2e^{-\mu_k}),\quad m=0,1,...,N. \nonumber
\end{align}
The value $v_k|_{\partial \Omega_k}$ is independent of $m$. In particular $\mu_{0,k}=\mu_k$. Thus the comparison of $\mu_{m,k}$ in (\ref{vk-bry}) proves
 (\ref{mu-lk}).
Next we observe that around $Q_l^k$
\begin{equation}\label{like-sing}
V_k(y)=\log \frac{e^{\mu_k}}{(1+\frac{|\tilde Q_l^k|^{2N}\mathfrak{h}_k(\delta_ke_1)e^{\mu_k}}
8 |y-\tilde Q_l^k|^2)^2}+\tilde\phi_l^k(y)+O(\mu_k^2e^{-\mu_k}),
\end{equation}
where $y\in B(e^{i\beta_l},\delta_0)$, $\tilde Q_l^k=e^{i\beta_l}+O(e^{-\mu_k})$,
\begin{equation}\label{har-a-l}
\tilde \phi_l^k(x)=\sum_{m=0,m\neq l}^N(-4)\log \frac{|y-e^{i\beta_m|}}{|e^{i\beta_m}-e^{i\beta_l}|},
\end{equation}
$\delta_0$ is a small positive number independent of $k$. The way to prove (\ref{like-sing}), by direct computation from the expression of $V_k$, is as follows: It is easy to see that
$$V_k(y)=-\frac 1{2\pi}\int_{\mathbb R^2}\log |y-\eta |\mathfrak{h}_k(\delta_ke_1)|\eta |^{2N}e^{V_k(\eta)}d\eta +C,\quad y\in \mathbb R^2. $$
Then $\mathfrak{h}_k(\delta_ke_1)e^{V_k}|y|^{2N}$ weakly converges to $8\pi \delta_{e^{i\beta_l}}$ in a small neighborhood of $e^{i\beta_l}$.  For $y\in \partial B(e^{i\beta_l},\delta)$ we have
$$V_k(y)=-\sum_{l=0}^N4\log |y-e^{i\beta_l}|+C_k+O(\mu_ke^{-\mu_k}). $$
From there we know that harmonic function around $e^{i\beta_l}$ that equals $0$ at $e^{i\beta_l}$ is $\phi_l^k$ in (\ref{har-a-l}). On the other hand the equation of $v_k$ around $e^{i\beta_l}$ is (\ref{eq-for-Vk}).
The standard expansion (see \cite{gluck}) for blowup solution leads to (\ref{like-sing}).

Since $Q_m^k-e^{i\beta_m}=E$, we can replace $\tilde \phi_l^k$ by $\phi_l^k$ and have
$$V_k(x)=\log \frac{e^{\mu_k}}{(1+\frac{\mathfrak{h}_k(\delta_ke_1)e^{\mu_k}}{8}|y-e^{i\beta_l}|^2)^2}+\phi_l^k+O(\mu_k^2e^{-\mu_k}) $$
 in $B(Q_l^k, e^{-\mu_k/2})$.
Thus in the region $B(Q_l^k, e^{-\mu_k/2})$, the comparison between $v_k$ and $V_k$ boils down to the evaluation of:

\begin{equation}\label{glob-a}
\log \frac{e^{\mu_{l,k}}}{(1+\frac{\mathfrak{h}_k(\delta_ke_1)e^{\mu_{l,k}}}8|y-e^{i\beta_l}-p_k|^2)^2}-\log \frac{e^{\mu_k}}{(1+\frac{\mathfrak{h}_k(\delta_ke_1)e^{\mu_k}}8|y-e^{i\beta_l}|^2)^2},
\end{equation}
for $|p_k|=E$. By elementary computation we see that the difference between the two terms in (\ref{glob-a}) is
$O(\mu_ke^{-\mu_k/2})$ if $|y-e^{i\beta_l}|\le C e^{-\mu_k/2}$.
On the other hand, for $Ce^{-\mu_k/2}<|y-e^{i\beta_l}|<\epsilon/2$, the comparison of expressions of $v_k$ and $U_k$ gives the difference upper bound as
$$O(e^{-\mu_k})|y-e^{i\beta_l}|^{-1}+O(\mu_k^2e^{-\mu_k}),\quad Ce^{-\mu_k/2}\le |y-e^{i\beta_l}|<\epsilon. $$
Moreover
\begin{equation}\label{around-l-vV}
v_k-V_k=O(\mu_k^2e^{-\mu_k})\quad \mbox{on}\quad \partial B(Q_l^k,\epsilon),\quad l=0,...,N.
\end{equation}
Also we observe from the expression of $V_k$ that
\begin{equation}\label{Vk-bry}
V_k(x)=-\mu_k-2\log \frac{\mathfrak{h}_k(\delta_ke_1)}8+4\log (N+1)-4(1+N)\log (\tau \delta_k^{-1})+E,
\end{equation}
for $ x\in \partial \Omega_k$,
thus
$$v_k-V_k=O(\mu_k^2e^{-\mu_k})\quad \mbox{on}\quad \partial \Omega_k. $$
Then
the closeness of $v_k$ and $V_k$ on $\Omega_k\setminus (\cup_lB(Q_l^k,\epsilon))$  can be obtained by a standard maximum principle argument: If we use $w_k$ to denote $v_k-V_k$:
$$w_k(z)=(v_k-V_k)(z), $$
then it is easy to see $w_k$ satisfies
$$|\Delta w_k(z)|\le Ce^{-\mu_k}|z|^{-4-2N},\quad \Omega_k\setminus B_2, \quad |w_k|\le C\mu_k^2e^{-\mu_k} \mbox{ on }\partial B_2\cup \partial \Omega_k, $$
 Then $|w_k|$ can be majorized by $Q(\mu_k^2e^{-\mu_k}-e^{-\mu_k}r^{-1-2N})$ for a large $Q>1$, which yields the smallness of $v_k-V_k$ on $\Omega_k\setminus B_2$ as a consequence.
Proposition \ref{vk-Vk} is established. $\Box$

\section{First crucial bound for $v_k-V_k$}

In this section we establish the first major estimate of $v_k-V_k$.
The main result in this section is

\begin{prop}\label{key-w8-8} Let $w_k=v_k-V_k$, then
$$|w_k(y)|\le C\delta_k, \quad y\in \Omega_k:=B(0,\tau \delta_k^{-1}). $$
\end{prop}

\noindent{\bf Proof of Proposition \ref{key-w8-8}:}

First we recall the equation for $v_k$ is (\ref{e-f-vk}),
$v_k=$ constant on $\partial B(0,\tau \delta_k^{-1})$. Moreover $v_k(e_1)=\mu_k$. Recall that $V_k$ defined in (\ref{def-Vk}) satisfies
$$\Delta V_k+\mathfrak{h}_k(\delta_ke_1)|y|^{2N}e^{V_k}=0,\quad \mbox{in}\quad \mathbb R^2, \quad \int_{\mathbb R^2}|y|^{2N}e^{V_k}<\infty, $$
$V_k$ has its local maximums at $e^{i\beta_l}$ for $l=0,...,N$ and $V_k(e_1)=\mu_k$.
For $|y|\sim \delta_k^{-1}$,
$$V_k(y)=-\mu_k-4(N+1)\log \delta_k^{-1}+C+O(\delta_k^{N+1})+O(e^{-\mu_k}). $$

Let $\Omega_k=B(0,\tau \delta_k^{-1})$,  we shall derive a precise, point-wise estimate of $w_k$ in $B_3\setminus \cup_{l=1}^{N}B(Q_l^k,\lambda)$ where $\lambda>0$ is a small number independent of $k$. Here we note that among $N+1$ local maximum points, we already have $e_1$ as a common local maximum point for both $v_k$ and $V_k$ and we shall prove that $w_k$ is very small in $B_3$ if we exclude all bubbling disks except the one around $e_1$. Before we carry out more specific computation we emphasize the importance of
\begin{equation}\label{control-e}
w_k(e_1)=|\nabla w_k(e_1)|=0.
\end{equation}
Now we write the equation of $w_k$ as
\begin{equation}\label{eq-wk}
\Delta w_k+\mathfrak{h}_k(\delta_k y)|y|^{2N}e^{\xi_k}w_k=(\mathfrak{h}_k(\delta_k e_1)-\mathfrak{h}_k(\delta_k y))|y|^{2N}e^{V_k}
\end{equation}
 in $\Omega_k$, where $\xi_k$ is obtained from the mean value theorem:
$$
e^{\xi_k(x)}=\left\{\begin{array}{ll}
\frac{e^{v_k(x)}-e^{V_k(x)}}{v_k(x)-V_k(x)},\quad \mbox{if}\quad v_k(x)\neq V_k(x),\\
\\
e^{V_k(x)},\quad \mbox{if}\quad v_k(x)=V_k(x).
\end{array}
\right.
$$
An equivalent form is
\begin{equation}\label{xi-k}
e^{\xi_k(x)}=\int_0^1\frac d{dt}e^{tv_k(x)+(1-t)V_k(x)}dt=e^{V_k(x)}\big (1+\frac 12w_k(x)+O(w_k(x)^2)\big ).
\end{equation}
For convenience we write the equation for $w_k$ as
\begin{equation}\label{eq-wk-2}
\Delta w_k+\mathfrak{h}_k(\delta_k y)|y|^{2N}e^{\xi_k}w_k=\delta_k\nabla \mathfrak{h}_k(\delta_k e_1)\cdot (e_1-y)|y|^{2N}e^{V_k}+E_1
\end{equation}
where $$E_1=O(\delta_k^2)|y-e_1|^2|y|^{2N}e^{V_k},\quad y\in \Omega_k. $$

Let $M_k=\max_{x\in \bar \Omega_k}|w_k(x)|$. We shall get a contradiction by assuming $M_k/\delta_k\to \infty$ at this moment.
Set
$$\tilde w_k(y)=w_k(y)/M_k,\quad x\in \Omega_k. $$
Clearly $\max_{x\in \Omega_k}|\tilde w_k(x)|=1$. The equation for $\tilde w_k$ is
\begin{equation}\label{t-wk}
\Delta \tilde w_k(y)+|y|^{2N}\mathfrak{h}_k(\delta_k e_1)e^{\xi_k}\tilde w_k(y)=\frac{\delta_k}{M_k}\nabla \mathfrak{h}_k(\delta_ke_1)\cdot (e_1-y)|y|^{2N}e^{V_k}+\tilde E_1,
\end{equation}
in $\Omega_k$,
where
\begin{equation}\label{t-ek}
\tilde E_1=o(\delta_k)|y-e_1|^2|y|^{2N}e^{V_k},\quad y\in \Omega_k.
\end{equation}

Now we give a more precise estimate of $e^{\xi_k}$. By Proposition \ref{vk-Vk}
\begin{equation}\label{xi-V-c}
\xi_k(y)=V_k(y)+\left\{\begin{array}{ll}
O(\mu_ke^{-\mu_k/2}),\quad |y-e_1|\le e^{-\mu_k/2},\\
\\
O(\mu_k^2e^{-\mu_k})|y-e_1|^{-1}, \,\, e^{-\mu_k/2}\le |y-e_1|\le \delta_0.
\end{array}
\right.
\end{equation}
Since $V_k$ is not exactly symmetric around $e_1$, we shall replace the re-scaled version of $V_k$ around $e_1$ by a radial function.
Let $U_k$ be solutions of
\begin{equation}\label{global-to-use}
\Delta U_k+\mathfrak{h}_k(\delta_ke_1)e^{U_k}=0,\quad \mbox{in}\quad \mathbb R^2, \quad U_k(0)=\max_{\mathbb R^2}U_k=0.
\end{equation}
Then we have
$$U_k(z)=\log \frac{1}{(1+\frac{\mathfrak{h}_k(\delta_ke_1)}{8}|z|^2)^2}$$
and
\begin{equation}\label{Vk-rad}
V_k(e_1+\epsilon_k z)+2\log \epsilon_k=U_k(z)+O(\epsilon_k)|z|+O(\mu_k^2\epsilon_k^2).
\end{equation}
Also we observe that
\begin{equation}\label{log-rad}
\log |e_1+\epsilon_ky|=O(\epsilon_k)|y|.
\end{equation}

Thus, the combination of (\ref{xi-V-c}), (\ref{Vk-rad}) and (\ref{log-rad}) gives
\begin{align}\label{xi-U}
&2N\log |e_1+\epsilon_kz|+\xi_k(e_1+\epsilon_k z)+2\log \epsilon_k-U_k(z)\\
=&O(\mu_k^2\epsilon_k)(1+|z|)\quad 0\le |z|<\delta_0 \epsilon_k^{-1}.
 \nonumber
\end{align}
Since we shall use the re-scaled version, based on (\ref{xi-U}) we have
\begin{equation}\label{xi-eU}
\epsilon_k^2 |e_1+\epsilon_k z|^{2N}e^{\xi_k(e_1+\epsilon_k z)}
= e^{U_k(z)}+O(\mu_k^2 \epsilon_k)(1+|z|)^{-3}
\end{equation}
Here we note that the estimate in (\ref{xi-U}) is not optimal.  In the following we shall put the proof of Proposition \ref{key-w8-8} into a few estimates. In the first estimate we prove

\begin{lem}\label{w-around-e1}
\begin{equation}\label{key-step-1}
\tilde w_k(y)=o(1),\quad \nabla \tilde w_k=o(1) \quad \mbox{in}\quad B(e_1,\delta)\setminus B(e_1,\delta/8)
\end{equation}
where $B(e_1,3\delta)$ does not include other blowup points.
\end{lem}

\noindent{\bf Proof of Lemma  \ref{w-around-e1}:}

If (\ref{key-step-1}) is not true, we have, without loss of generality that $\tilde w_k\to c>0$. Note that $\tilde w_k$ tends to a global harmonic function with removable singularity. So $\tilde w_k$ tends to constant. Here we assume $c>0$ but the argument for $c<0$ is the same. Let
\begin{equation}\label{w-ar-e1}
W_k(z)=\tilde w_k(e_1+\epsilon_kz), \quad \epsilon_k=e^{-\frac 12 \mu_k},
\end{equation}
then if we use $W$ to denote the limit of $W_k$, we have
$$\Delta W+e^UW=0, \quad \mathbb R^2, \quad |W|\le 1, $$
and $U$ is a solution of $\Delta U+e^U=0$ in $\mathbb R^2$ with $\int_{\mathbb R^2}e^U<\infty$. Since $0$ is the local maximum of $U$,
$$U(z)=\log \frac{1}{(1+\frac 18|z|^2)^2}. $$
Here we further claim that $W\equiv 0$ in $\mathbb R^2$ because $W(0)=|\nabla W(0)|=0$, a fact well known based on the classification of the kernel of the linearized operator. Going back to $W_k$, we have
$$W_k(z)=o(1),\quad |z|\le R_k \mbox{ for some } \quad R_k\to \infty. $$

Based on the expression of $\tilde w_k$, (\ref{Vk-rad}) and (\ref{xi-eU}) we write the equation of $W_k$ as
\begin{equation}\label{e-Wk}
\Delta W_k(z)+\mathfrak{h}_k(\delta_ke_1)e^{U_k(z)}W_k(z)=-\frac{\delta_k}{M_k}\nabla\mathfrak{h}_k(\delta_ke_1)\cdot z\epsilon_ke^{U_k(z)}+E_2^k,
\end{equation}
for $|z|<\delta_0 \epsilon_k^{-1}$ where
\begin{equation*}
E_2^k(z)=o(1)\mu_k^2\epsilon_k(1+|z|)^{-3}.
\end{equation*}

Let
\begin{equation}\label{for-g0}
g_0^k(r)=\frac 1{2\pi}\int_0^{2\pi}W_k(r,\theta)d\theta.
\end{equation}
Then clearly $g_0^k(r)\to c>0$ for $r\sim \epsilon_k^{-1}$.
 The equation for $g_0^k$ is
\begin{align*}
&\frac{d^2}{dr^2}g_0^k(r)+\frac 1r \frac{d}{dr}g_0^k(r)+\mathfrak{h}_k(\delta_ke_1)e^{U_k(r)}g_0^k(r)=\tilde E_0^k(r)\\
&g_0^k(0)=\frac{d}{dr}g_0^k(0)=0.
\end{align*}
where $\tilde E_0^k(r)$ has the same upper bound as that of $E_2^k(r)$:
$$|\tilde E_0^k(r)|\le C\mu_k^2\epsilon_k(1+r)^{-3}. $$

For the homogeneous equation, the two fundamental solutions are known: $g_{01}$, $g_{02}$, where
$$g_{01}=\frac{1-c_1r^2}{1+c_1r^2},\quad c_1=\frac{\mathfrak{h}_k(\delta_ke_1)}8.$$
By the standard reduction of order process, $g_{02}(r)=O(\log r)$ for $r>1$.
Then it is easy to obtain, assuming $|W_k(z)|\le 1$, that
\begin{align*}
|g_0(r)|\le C|g_{01}(r)|\int_0^r s|\tilde E_0^k(s) g_{02}(s)|ds+C|g_{02}(r)|\int_0^r s|g_{01}(s)\tilde E_0^k(s)|ds\\
\le C\mu_k^2\epsilon_k\log (2+r). \quad 0<r<\delta_0 \epsilon_k^{-1}.
\end{align*}
Clearly this is a contradiction to (\ref{for-g0}). We have proved $c=0$, which means $\tilde w_k=o(1)$ in $B(e_1, \delta_0)\setminus B(e_1, \delta_0/8)$.
Then it is easy to use the equation for $\tilde w_k$ and standard Harnack inequality to prove
$\nabla \tilde w_k=o(1)$ in the same region.
Lemma \ref{w-around-e1} is established. $\Box$

\medskip

The second estimate is a more precise description of $\tilde w_k$ around $e_1$:
\begin{lem}\label{t-w-1-better} For any given $\sigma\in (0,1)$ there exists $C>0$ such that
\begin{equation}\label{for-lambda-k}
|\tilde w_k(e_1+\epsilon_kz)|\le C\epsilon_k^{\sigma} (1+|z|)^{\sigma},\quad 0<|z|<\tau \epsilon_k^{-1}.
\end{equation}
for some $\tau>0$.
\end{lem}

\noindent{\bf Proof of Lemma \ref{t-w-1-better}:} Let $W_k$ be defined as in (\ref{w-ar-e1}) we recall the equation for $W_k$ as (\ref{e-Wk}).
Here we let $\Omega_{wk}=B(0,\tau \epsilon_k^{-1})$ where $\tau\in (0,1)$ is chosen such that $W_k=o(1)$ on $\partial \Omega_{wk}$.
  and
$$\Lambda_k=\max_{z\in \Omega_{wk}}\frac{|W_k(z)|}{\epsilon_k^{\sigma}(1+|z|)^{\sigma}}. $$
If (\ref{for-lambda-k}) does not hold, $\Lambda_k\to \infty$ and we use $z_k$ to denote where $\Lambda_k$ is attained. Note that because of the smallness of $W_k$ on $\partial \Omega_{wk}$, $z_k$ is an interior point. Let
$$g_k(z)=\frac{W_k(z)}{\Lambda_k (1+|z_k|)^{\sigma}\epsilon_k^{\sigma}},\quad z\in \Omega_{k,1}, $$
we see immediately that
\begin{equation}\label{g-sub-linear}
|g_k(z)|=\frac{|W_k(z)|}{\epsilon_k^{\sigma}\Lambda_k(1+|z|)^{\sigma}}\cdot \frac{(1+|z|)^{\sigma}}{(1+|z_k|)^{\sigma}}\le  \frac{(1+|z|)^{\sigma}}{(1+|z_k|)^{\sigma}}.
\end{equation}
Note that $\sigma$ can be as close to $1$ as needed. The equation of $g_k$ is
$$\Delta g_k(z)+\mathfrak{h}_k(\delta_k e_1)e^{\xi_k}g_k=o(\epsilon_k^{1-\sigma})\frac{(1+|z|)^{-3}}{(1+|z_k|)^{\sigma}}, \quad \mbox{in}\quad \Omega_{k,1}. $$
Then we can obtain a contradiction to $|g_k(z_k)|=1$ as follows: If $\lim_{k\to \infty}z_k=P\in \mathbb R^2$, this is not possible because that fact that $g_k(0)=|\nabla g_k(0)|=0$ and the sub-linear growth of $g_k$ in (\ref{g-sub-linear}) implies that $g_k\to 0$ over any compact subset of $\mathbb R^2$. So we have $|z_k|\to \infty$. But this would lead to a contradiction again by using the Green's representation of $g_k$:
\begin{align} \label{temp-1}
&\pm 1=g_k(z_k)=g_k(z_k)-g_k(0)\\
&=\int_{\Omega_{k,1}}(G_k(z_k,\eta)-G_k(0,\eta))(\mathfrak{h}_k(\delta_k e_1)e^{\xi_k}g_k(\eta)+o(\epsilon_k^{1-\sigma})\frac{(1+|\eta |)^{-3}}{(1+|z_k|)^{\sigma}})d\eta+o(1).\nonumber
\end{align}
where $G_k(y,\eta)$ is the Green's function on $\Omega_{k,1}$ and $o(1)$ in the equation above comes from the smallness of $W_k$ on $\partial \Omega_{k,1}$. Let $L_k=\tau\epsilon_k^{-1}$. From (\ref{expre-G}) we have
$$G_k(z_k,\eta)-G_k(0,\eta)=-\frac{1}{2\pi}\log |z_k-\eta |+\frac 1{2\pi}\log |\frac{z_k}{|z_k|}-\frac{\eta z_k}{L_k}|+\frac 1{2\pi}\log |\eta |. $$
Using this expression in (\ref{temp-1}) we obtain from elementary computation that the right hand side of (\ref{temp-1}) is $o(1)$, a contradiction to $|g_k(z_k)|=1$. Lemma \ref{t-w-1-better} is
established. $\Box$

\medskip

The smallness of $\tilde w_k$ around $e_1$ can be used to obtain the following third key estimate:
\begin{lem}\label{small-other}
\begin{equation}\label{key-step-2}
\tilde w_k=o(1)\quad \mbox{in}\quad B(e^{i\beta_l},\tau)\quad l=1,..,N.
\end{equation}
\end{lem}

\noindent{\bf Proof of Lemma \ref{small-other}:}
We abuse the notation $W_k$ by defining it as
$$W_k(z)=\tilde w_k(e^{i\beta_l}+\epsilon_k z),\quad z\in \Omega_{k,l}:=B(0,\tau \epsilon_k^{-1}). $$
Here we point out that the smallness of $\delta_k$ (see \ref{small-delta})implies that $\epsilon_k^{-1}|Q_l^k-e^{i\beta_l}|\to 0$. So the scaling around $e^{i\beta_l}$ or $Q_l^k$ does not affect the
limit function.
$$|e^{i\beta_l}+\epsilon_kz|^{2N}\mathfrak{h}_k(\delta_ke_1)e^{\xi_k(e^{i\beta_l}+\epsilon_kz)}\to e^{U(z)} $$
where $U(z)$ is a solution of
$$\Delta U+e^U=0,\quad \mbox{in}\quad \mathbb R^2, \quad \int_{\mathbb R^2}e^U<\infty. $$
Here we recall that $\lim_{k\to \infty} \mathfrak{h}_k(\delta_k e_1)=1$.
Since $W_k$ converges to a solution of the linearized equation:
$$\Delta W+e^UW=0, \quad \mbox{in}\quad \mathbb R^2. $$
$W$ can be written as a linear combination of three functions:
$$W(x)=c_0\phi_0+c_1\phi_1+c_2\phi_2, $$
where
$$\phi_0=\frac{1-\frac 18 |x|^2}{1+\frac 18 |x|^2} $$
$$\phi_1=\frac{x_1}{1+\frac 18 |x|^2},\quad \phi_2=\frac{x_2}{1+\frac 18|x|^2}. $$
The remaining part of the proof consisting of proving $c_0=0$ and $c_1=c_2=0$. First we prove $c_0=0$.

\noindent{\bf Step one: $c_0=0$.}
First we write the equation for $W_k$ in a convenient form. Since
$$|e^{i\beta_l}+\epsilon_kz|^{2N}\mathfrak{h}_k(\delta_ke_1)=\mathfrak{h}_k(\delta_ke_1)+O(\epsilon_k z),$$
and
$$\epsilon_k^2e^{\xi_k(e^{i\beta_l}+\epsilon_kz)}=e^{U_k(z)}+O(\epsilon_k)(1+|z|)^{-3}. $$
Based on (\ref{t-wk}) we write the equation for $W_k$ as
\begin{equation}\label{around-l}
\Delta W_k(z)+\mathfrak{h}_k(\delta_ke_1)e^{U_k}W_k=E_l^k(z)
\end{equation}
where
$$E_l^k(z)=O(\epsilon_k^{1/2})(1+|z|)^{-3}+\sigma_k\mathfrak{h}_k(\delta_ke_1)e^{U_k}\quad \mbox{in}\quad \Omega_{k,l} $$
where $\sigma_k=\frac{\delta_k}{M_k}\to 0$.
In order to prove $c_0=0$, the key is to control the derivative of $W_0^k(r)$ where
$$W_0^k(r)=\frac 1{2\pi r}\int_{\partial B_r} W_k(re^{i\theta})dS, \quad 0<r<\tau \epsilon_k^{-1}. $$
To obtain a control of $\frac{d}{dr}W_0^k(r)$ we use $\phi_0^k(r)$ as the radial solution of
$$\Delta \phi_0^k+\mathfrak{h}_k(\delta_k e_1)e^{U_k}\phi_0^k=0, \quad \mbox{in }\quad \mathbb R^2. $$
When $k\to \infty$, $\phi_0^k\to c_0\phi_0$. Thus using the equation for $\phi_0^k$ and $W_k$, we have
\begin{equation}\label{c-0-pf}
\int_{\partial B_r}(\partial_{\nu}W_k\phi_0^k-\partial_{\nu}\phi_0^kW_k)=o(\epsilon_k^{1/2})+\sigma_k\int_{B_r}\mathfrak{h}_k(\delta_k e_1)e^{U_k}\phi_k. \end{equation}
The last term above is
$$\sigma_k\int_{B_r}\mathfrak{h}_k(\delta_ke_1)e^{U_k}\phi_k=\sigma_k\int_{\partial B_r}\partial_{\nu}\phi_0^k=o(1)r^{-3},\quad r>1. $$
Thus from (\ref{c-0-pf}) we have
\begin{equation}\label{W-0-d}
\frac{d}{dr}W_0^k(r)=o(\epsilon_k^{1/2})/r+O(1/r^3),\quad 1<r<\tau \epsilon_k^{-1}.
\end{equation}
Since we have known that
$$W_0^k(\tau \epsilon_k^{-1})=o(1). $$
By the fundamental theorem of calculus we have
$$W_0^k(r)=W_0^k(\tau\epsilon_k^{-1})+\int_{\tau \epsilon_k^{-1}}^r(\frac{o(\epsilon_k^{\frac 12})}{s}+O(s^{-3}))ds=O(1/r^2), \quad r\sim 1,\quad r \mbox{ is large}. $$
$c_0=0$ is established because $W_0^k(r)\to c_0\phi_0$.

\medskip

\noindent{\bf Step two $c_1=c_2=0$}.
 We first observe that Lemma \ref{small-other} follows from this. Indeed, once we have proved $c_1=c_2=c_0=0$ around each $e^{i\beta_l}$, it is easy to use maximum principle to prove $\tilde w_k=o(1)$ in $B_3$ using $\tilde w_k=o(1)$ on $\partial B_3$ and the Green's representation of $\tilde w_k$. The smallness of $\tilde w_k$ immediately implies $\tilde w_k=o(1)$ in $B_R$ for any fixed $R>>1$. Outside $B_R$, a crude estimate of $v_k$ is
  $$v_k(y)\le -\mu_k-4(N+1)\log |y|+C, \quad 3<|y|<\tau \delta_k^{-1}. $$
  Using this and the Green's representation of $w_k$ we can first observe that the oscillation on each $\partial B_r$ is $o(1)$ ($R<r<\tau \delta_k^{-1}/2$) and then by the Green's representation of $\tilde w_k$ and fast decay rate of $e^{V_k}$ we obtain $\tilde w_k=o(1)$ in $B(0,\tau \delta_k^{-1})$. A contradiction to $\max |\tilde w_k|=1$.

 There are $N+1$ local maximums with one of them being $e_1$. Correspondingly there are $N+1$ global solutions $V_{l,k}$ that
approximate $v_k$ accurately near $Q_l^k$ for $l=0,...,N$. Note that $Q_0^k=e_1$. For $V_{l,k}$ the expression is
$$V_{l,k}=\log \frac{e^{\mu_l^k}}{(1+\frac{e^{\mu_l^k}}{D_l^k}|y^{N+1}-(e_1+p_l^k)|^2)^2},\quad l=0,...,N, $$
where $p_l^k=E$ and $D_l^k=8(N+1)^2/\mathfrak{h}_k(\delta_kQ_l^k)$. Note that a simple computation gives
$$Q_l^k=e^{i\beta_l}+\frac{p_l^ke^{i\beta_l}}{N+1}+O(|p_l^k|^2),\quad \beta_l=\frac{2l\pi}{N+2}.$$
Let $M_{l,k}$ be the maximum of $|v_k-V_{l,k}|$ and we claim that all these $M_{l,k}$ are comparable:
\begin{equation}\label{M-comp}
M_{l,k}\sim M_{s,k},\quad \forall s\neq l.
\end{equation}
The proof of (\ref{M-comp}) is as follows: We use $L_{s,l}$ to denote the limit of $v_k-V_{l,k}$ around $Q_s^k$:
$$\frac{(v_k-V_{l,k})(Q_s^k+\epsilon_kz)}{M_{l,k}}=L_{s,l}+o(1),\quad |z|\le \tau \epsilon_k^{-1} $$
where
$$ L_{s,l}=c_{1,s,l}\frac{z_1}{1+\frac 18 |z|^2}+c_{2,s,l}\frac{z_2}{1+\frac 18 |z|^2},\quad \mbox{and}\quad L_{l,l}=0, \quad s=0,...,N. $$
If all $c_{1,s,l}$ and $c_{2,s,l}$ are zero there is nothing to prove. So at least one of them is not zero.
For each $s\neq l$, by Lemma \ref{t-w-1-better} we have
\begin{equation}\label{Q-bad}
v_k(Q_s^k+\epsilon_kz)-V_{s,k}(Q_s^k+\epsilon_kz)=O(\epsilon_k^{\sigma})(1+|z|)^{\sigma} M_{s,k},\quad |z|<\tau \epsilon_k^{-1}.
\end{equation}
Let $M_k=\max_{i}M_{i,k}$ ($i=0,...,N$) and we suppose $M_k=M_{l,k}$. Then to determine $L_{s,l}$ we see that
$$\frac{v_k(Q_s^k+\epsilon_k z)-V_{l,k}(Q_s^k+\epsilon_kz)}{M_k}=o(\epsilon_k^{\sigma})(1+|z|)^{\sigma}+\frac{V_{s,k}(Q_s^k+\epsilon_kz)-V_{l,k}(Q_s^k+\epsilon_kz)}{M_k}. $$
We write $V_{s,k}(y)-V_{l,k}(y)$ as
$$V_{s,k}(y)-V_{l,k}(y)=\mu_{s,k}-\mu_{l,k}+2A-A^2+O(|A|^3) $$
where
$$A(y)=\frac{\frac{e^{\mu_l^k}}{D_l^k}|y^{N+1}-e_1-p_l^k|^2-\frac{e^{\mu_s^k}}{D_s^k}|y^{N+1}-e_1-p_s^k|^2}{1+\frac{e^{\mu_s^k}}{D_s^k}|y^{N+1}-e_1-p_s^k|^2}.$$
From $A$ we have
\begin{align}\label{late-1}
&V_{s,k}(Q_s^k+\epsilon_kz)-V_{l,k}(Q_s^k+\epsilon_kz)\\
=&\phi_1+\phi_2+\phi_3+\phi_4+\mathfrak{R},\nonumber
\end{align}
where
\begin{align*}
&\phi_1=(\mu_{l,k}-\mu_{s,k})(1-\frac{(N+1)^2}{D_s^k}|z+\frac{N}2z^2e^{-i\beta_s}\epsilon_k|^2)/B, \\
&\phi_2=\frac {2(N+1)^2}{D_s^k} \delta_k\nabla\mathfrak{h}_k(\delta_k Q_s^k)(Q_l^k-Q_s^k)|z|^2/B\\
&\phi_3=\frac{4e^{\mu_{l,k}-\mu_{s,k}}(N+1)^2}{D_s^k B}Re((z+\frac N2\epsilon_ke^{-i\beta_s}z^2)(\frac{\bar p_s^k-\bar p_l^k}{\epsilon_k}e^{-i\beta_s}))
\\
&\phi_4=\frac{|p_s^k-p_l^k|^2}{\epsilon_k^2}\bigg (\frac{2}{D_s^k  B}-\frac{2(N+1)^2|z|^2}{D_s^2 B^2}
-\frac{2(N+1)^2}{D_s^2B}|z|^2\cos(2\theta-2\theta_{st}-2\beta_s)\bigg ),\\
&B=1+\frac{(N+1)^2}{D_s^k}|z+\frac N2z^2e^{-i\beta_s}
\epsilon_k|^2,
\end{align*}
and $\mathfrak{R}_k$ is the collections of other insignificant terms. Here we note that if we set $\epsilon_{l,k}=e^{-\mu_l^k/2}$, there is no essential difference between $\epsilon_{l,k}$ and $\epsilon_k=e^{-\frac 12\mu_{1,k}}$ because $\epsilon_{l,k}=\epsilon_k+O(\epsilon_k E)$. If $|\mu_{s,k}-\mu_{l,k}|/M_k\ge C$ there is no way to obtain a limit in the form of $L_{s,l}$ mentioned before. Thus we must have $|\mu_{s,k}-\mu_{l,k}|/M_k\to 0$. After simplification (see $\phi_3$ of (\ref{late-1})) we have
\begin{equation}\label{c-12}
c_{1,l,s}=\lim_{k\to \infty}\frac{|p_s^k-p_l^k|}{2(N+1)M_k\epsilon_k}\cos(\beta_s+\theta_{sl}),\quad c_{2,l,s}=\lim_{k\to \infty}
\frac{|p_s^k-p_l^k|}{2(N+1)\epsilon_k M_k}sin(\beta_s+\theta_{sl})
\end{equation}
and $\theta_{sl}$ comes from $p_s^k-p_l^k=|p_s^k-p_l^k|e^{i\theta_{sl}}$. We omit $k$ for convenience.
It is also important to observe that even if $M_k=o(\epsilon_k)$ we still have $M_k\sim \max_{s}|p_s^k-p_l^k|/\epsilon_k$. Since each $|p_l^k|=E$, an upper bound for $M_k$ is $M_k\le C\mu_k\epsilon_k$.

Equation (\ref{c-12}) gives us the important observation: $|c_{1,l,s}|+|c_{2,l,s}|\sim |p_s^k-p_l^k|/M_k$. So whenever $|c_{1,l,s}|+|c_{2,l,s}|\neq 0$ we have
$|p_s^k-p_l^k|\sim \epsilon_k M_k$. In other words for each $l$, $\epsilon_k M_{l,k}\sim \max_{t\neq l}|p_t^k-p_l^k|$.  Hence for any $t$, if $|p_t^k-p_l^k|\sim \epsilon_k M_k$, let $M_{t,k}$ be the maximum of $|v_k-V_{t,k}|$, we have $M_{t,k}\sim M_k$. If all $|p_t^k-p_l^k|\sim \epsilon_k M_k$ (\ref{M-comp}) is proved. So we assume that at least one of them is too close to $p_l^k$, say $|p_s^k-p_l^k|=o(1)\epsilon_k M_k$. We choose $t\neq l$ such that $|p_t^k-p_l^k|\sim \epsilon_k M_k$. Then
$$|p_s^k-p_t^k|\ge |p_t^k-p_l^k|-|p_s^k-p_l^k|\ge \frac 12|p_t^k-p_l^k|\sim \epsilon_k M_k. $$
This means $M_{s,k}\sim M_k$. Thus (\ref{M-comp}) is established.  From now on for convenience we shall just use $M_k$.

Set $w_{l,k}=(v_k-V_{l,k})/M_{k}$, then we have $w_{l,k}(Q_l^k)=|\nabla w_{l,k}(Q_l^k)|=0$,
\begin{equation}\label{aroud-s-1}
\lim_{k\to \infty}w_{l,k}(Q_s^k+\epsilon_k z)/M_k=\frac{c_{1,l,s}z_1+c_{2,l,s}z_2}{1+\frac 18 |z|^2}
\end{equation}
and around each $Q_s^k$ (\ref{Q-bad}) holds with $M_{s,k}$ replaced by $M_k$.

Now we use the estimates above to evaluate $\nabla w_{l,k}(Q_l^k)=0$. By elementary computation,
\begin{align}\label{key-estimate}
&\nabla w_{l,k}(Q_l^k)\\
=&-\frac 1{2\pi}\int_{\Omega_k}\frac{Q_l^k-\eta}{|Q_l^k-\eta |^2}(w_{l,k}(\eta)\mathfrak{h}_k(\delta_kQ_l^k)+\sigma_{k}\nabla \mathfrak{h}_k(\delta_k Q_l^k)(\eta -Q_l^k))|\eta |^{2N}e^{V_k}d\eta +o(\epsilon_k).\nonumber
\end{align}
Since the left hand side is $0$, the right hand side gives
\begin{align}\label{eq-ql}
\sum_{s\neq l}\bigg  (\int_{B(Q_s^k,\tau)}\frac{Q_l^k-\eta}{|Q_l^k-\eta |^2}w_{l,k}(\eta)\mathfrak{h}_k(\delta_kQ_l^k)|\eta |^{2N}e^{V_k}d\eta\\
+8\pi\sigma_{k}\nabla \log \mathfrak{h}_k(\delta_k Q_s^k)
(Q_s^k-Q_l^k)\bigg )=o(\epsilon_k), \nonumber
\end{align}
where $\sigma_{k}=\delta_k/M_{k}\to 0$. To evaluate the first term in (\ref{eq-ql}), around each $Q_s^k$ we write
$$\frac{Q_l^k-\eta}{|Q_l^k-\eta |^2}=\frac{Q_l^k-Q_s^k}{|Q_l^k-Q_s^k|^2}+F_{ls}(\eta) $$
where
$$F_{ls}(\eta)=(F_{ls}^1(\eta),F_{ls}^2(\eta))= \frac{Q_l^k-\eta}{|Q_l^k-\eta |^2}-\frac{Q_l^k-Q_s^k}{|Q_l^k-Q_s^k|^2}.$$
We use (\ref{late-1}) to obtain
\begin{align}\label{late-2-r}
\int_{B(Q_s^k,\tau)}(w_{l,k}(\eta)\mathfrak{h}_k(\delta_k Q_l^k)|\eta |^{2N}e^{V_k}d\eta\\
+8\pi\sigma_{k}\nabla \log \mathfrak{h}_k(\delta_k Q_s^k)
(Q_s^k-Q_l^k)=O(\epsilon_k^{\sigma}) \nonumber
\end{align}

Note that in the computation above, the terms of $\phi_1$ and $\phi_3$ lead to $o(\epsilon_k)$, the integration involving $\phi_2$ cancels with the second term of (\ref{late-2-r}). The integration involving $\phi_4$ provides the leading term. More detailed information is the following:
First for a global solution
$$V_{\mu,p}=\log \frac{e^{\mu}}{(1+\frac{e^{\mu}}{\lambda}|z^{N+1}-p|^2)^2}$$ of
$$\Delta V_{\mu,p}+\frac{8(N+1)^2}{\lambda}|z|^{2N}e^{V_{\mu,p}}=0,\quad \mbox{in }\quad \mathbb R^2, $$ we have
\begin{equation}\label{inte-eq-1}
\int_{\mathbb R^2}\partial_{\mu}V_{\mu,p}|z|^{2N}e^{V_{\mu,p}}=\int_{\mathbb R^2}\frac{(1-\frac{e^{\mu}}{\lambda}|z^{N+1}-P|^2)|z|^{2N}}{(1+\frac{e^{\mu}}{\lambda}|z^{N+1}-P|^2)^3}dz=0.
\end{equation}
From this we use standard estimate that includes some cancellation to have
$$\int_{B(0,\tau\epsilon_k^{-1})}\frac{\phi_1}{M_k}B^{-3}=o(\epsilon_k).$$
 From $V_{\mu,p}$  we also have
$$\int_{\mathbb R^2}\partial_{P}V_{\mu,p}|y|^{2N}e^{V_{\mu,p}}=\int_{\mathbb R^2}\partial_{\bar P}V_{\mu,p}|y|^{2N}e^{V_{\mu,p}}=0, $$
which gives
\begin{equation}\label{inte-eq-2}
\int_{\mathbb R^2}\frac{\frac{e^{\mu}}{\lambda}(\bar z^{N+1}-\bar P)|z|^{2N}}{(1+\frac{e^{\mu}}{\lambda}|z^{N+1}-P|^2)^3}=\int_{\mathbb R^2}\frac{\frac{e^{\mu}}{\lambda}( z^{N+1}- P)|z|^{2N}}{(1+\frac{e^{\mu}}{\lambda}|z^{N+1}-P|^2)^3}=0.
\end{equation}
From here we use scaling and cancellation to have
$$\int_{B(0,\tau\epsilon_k^{-1})}\frac{\phi_3}{M_k}B^{-3}=o(\epsilon_k).$$ Thus (\ref{late-2-r}) holds.
Equation (\ref{late-2-r}) also leads to a more accurate estimate of $w_{l,k}$ in regions between bubbling disks. By the Green's representation formula of $w_{l,k}$ it is easy to have
\begin{align*}
& w_{l,k}(y)\\
=&-\frac{1}{2\pi}\int_{\Omega_k}\log \frac{|y-\eta |}{|Q_l^k-\eta |}( w_{l,k}(\eta )\mathfrak{h}_k(\delta_k Q_l^k)+\nabla \mathfrak{h}_k(\delta_k Q_l^k)(\eta-Q_l^k)
|\eta |^{2N}e^{V_k}d\eta+o(\epsilon_k)
\end{align*}
Writing
$$\log \frac{|y-\eta |}{|Q_l^k-\eta |}=\log \frac{|y-Q_s^k|}{|Q_l^k-Q_s^k|}+(\log \frac{|y-\eta |}{|Q_l^k-\eta |}-\log \frac{|y-Q_s^k|}{|Q_l^k-Q_s^k|}),$$
the integration related to the second term is $O(\epsilon_k)$. The integration involving the first term is $O(\epsilon_k^{\sigma})$ by (\ref{late-2-r}).
Therefore
$$|w_{l,k}(y)|=o(\epsilon_k^{\sigma}),\quad y\in B_3\setminus \cup_{s=0}^NB(Q_s^k,\tau) $$
for $\sigma\in (0,1)$. Thus this extra control of $w_{l,k}$ away from bubbling disks gives a better estimate than (\ref{Q-bad}) around each $Q_s^k$:
\begin{equation}\label{much-better-s}
|w_{s,k}(Q_s^k+\epsilon_k z)|/M_k\le o(\epsilon_k)(1+|z|)^{\sigma},\quad |z|<\tau \epsilon_k^{-1}.
\end{equation}

From the decomposition in (\ref{late-1}) we can now compute (\ref{key-estimate}) in more detail:

\begin{align}\label{late-2}
\int_{B(Q_s^k,\tau)}w_{l,k}(\eta)\mathfrak{h}_k(\delta_k Q_l^k)|\eta |^{2N}e^{V_k}d\eta+8\pi\sigma_{k}\nabla \log \mathfrak{h}_k(\delta_k Q_s^k)
(Q_s^k-Q_l^k)\\
=\frac{\pi}{(N+1)^2}\frac{|p_s^k-p_l^k|^2}{\epsilon_k^2 M_{k}^2}M_k+O(\epsilon_k^{\sigma}).\nonumber
\end{align}

\begin{align}\label{key-est-2}
&\int_{B(Q_s^k,\tau)}F_{ls}^1(\eta)(w_{l,k}(\eta)\mathfrak{h}_k(\delta_k Q_l^k)|\eta |^{2N}e^{V_k}d\eta \\
&=\epsilon_k\int_{B(0,\tau\epsilon_k^{-1})}(\partial_1 F_{ls}^1(Q_s^k)z_1+\partial_2F^1_{ls}(Q_s^k)z_2)\frac{w_{l,k}(Q_s^k+\epsilon_kz)}{(1+\frac{1}8|z|^2)^2}dz+o(\epsilon_k).\nonumber \\
&=4\pi\epsilon_k \frac{-c_{1,l,s}\cos(\beta_s+\beta_l)-c_{2,l,s}\sin(\beta_l+\beta_s)}{\sin^2(\frac{\beta_l-\beta_s}2)}+o(\epsilon_k).\nonumber
\end{align}

\begin{align}\label{key-est-3}
&\int_{B(Q_s^k,\tau)}F_{ls}^2(\eta)(w_{l,k}(\eta)\mathfrak{h}_k(\delta_k Q_l^k)|\eta |^{2N}e^{V_k}d\eta \\
&=\epsilon_k\int_{B(0,\tau\epsilon_k^{-1})}(\partial_1 F_{ls}^2(Q_s^k)z_1+\partial_2F^2_{ls}(Q_s^k)z_2)\frac{w_{l,k}(Q_s^k+\epsilon_kz)}{(1+\frac{1}8|z|^2)^{2}}dz+o(\epsilon_k).\nonumber \\
&=4\pi\epsilon_k \frac{-c_{1,l,s}\sin(\beta_s+\beta_l)+c_{2,l,s}\cos(\beta_l+\beta_s)}{\sin^2(\frac{\beta_l-\beta_s}2)}+o(\epsilon_k).\nonumber
\end{align}

Using (\ref{late-2}),(\ref{key-est-2}),(\ref{key-est-3}) in (\ref{eq-ql}), and by the expressions of $c_{1,l,s},c_{2,l,s}$ in (\ref{c-12}) we have
\begin{equation*}
\sum_{s\neq l}\bigg (\frac{Q_l^k-Q_s^k}{|Q_l^k-Q_s^k|^2}\frac{\pi}{(N+1)^2}\frac{|p_s^k-p_l^k|^2}{\epsilon_k^2M_k^2}M_k-
\frac{2\pi\epsilon_k}{\sin^2(\frac{\beta_l-\beta_s}{2})}e^{i(\beta_l-\theta_{sl})}\frac{|p_s^k-p_l^k|}{\epsilon_k M_k}\bigg)=o(\epsilon_k).
\end{equation*}

Since $Q_s^k=e^{i\beta_s}+E$, the equation above can be simplified to
\begin{align}\label{key-6}
\sum_{s\neq l}\bigg ((\frac 12-\frac{i\sin(\beta_l-\beta_s)}{4\sin^2(\frac{\beta_s-\beta_l}2)})\frac{\pi}{(N+1)^2}(\frac{|p_l^k-p_s^k|}{\epsilon_k M_k})^2M_k\\
-\frac{2\pi\epsilon_ke^{-i\theta_{sl}}}{\sin^2(\frac{p_s^k-p_l^k}2)}\frac{|p_s^k-p_l^k|}{\epsilon_k M_k}\bigg )=o(\epsilon_k). \nonumber
\end{align}
Note that $e^{-i\theta_{sl}}|p_s^k-p_l^k|=\bar p_s^k-\bar p_l^k$. Taking the sum of all $l$ in (\ref{key-6}) we have
$$CM_k\sum_{s\neq l}(|p_s^k-p_l^k|^2/(\epsilon_k M_k)^2)=o(\epsilon_k), \quad C>0. $$
If $M_k\ge C\epsilon_k$ for some $C>0$, we have $|c_{1,s,l}|+|c_{2,s,l}|=0$ for all $s,l$. Finally if $M_k=o(\epsilon_k)$, the first term in (\ref{key-6}) is ignored and (\ref{key-6}) becomes
$$\sum_{s\neq l}\frac{\bar p_s^k-\bar p_l^k}{\sin^2(\frac{\beta_s-\beta_l}2) \epsilon_k M_k}=o(1),\quad l=0,1,...,N. $$
From here it is also easy to obtain all $|p_l^k|/(\epsilon_k M_k)=o(1)$, which implies $c_{1,s,l}=c_{2,s,l}=0$ for all $s, l$. Thus
Lemma \ref{small-other} is established. $\Box$

\medskip

Proposition \ref{key-w8-8} is an immediate consequence of Lemma \ref{small-other}.  $\Box$.
\section{Proof of Theorem \ref{main-thm}}

Let $\hat w_k=w_k/\delta_k$. Then the equation for $\hat w_k$ is
\begin{equation}\label{hat-w}
\Delta \hat w_k+|y|^{2N}e^{\xi_k}\hat w_k=\nabla\mathfrak{h}_k(0)\cdot (e_1-y)|y|^{2N}e^{V_k}+O(\delta_k)e^{V_k}|y-e_1|^{2},
\end{equation}
in $\Omega_k$.
By Proposition \ref{key-w8-8}, $|\hat w_k(y)|\le C$. Before we carry out the remaining part of the proof we observe that $\hat w_k$ converges to a harmonic function in $\mathbb R^2$ minus finite singular points. Since $\hat w_k$ is bounded, all these singularities are removable. Thus $\hat w_k$ converges to a constant. Based on the information around $e_1$, we shall prove that this constant is $0$. However, looking at the right hand side the equation,
$$\nabla\mathfrak{h}_k(0)\cdot (e_1-y)|y|^{2N}e^{V_k}\rightharpoonup \sum_{l=1}^N8\pi \nabla \mathfrak{h}_k(0)\cdot (e_1-e^{i\beta_l})\delta_{e^{i\beta_l}}. $$
If $\nabla \mathfrak{h}_k(0)\neq 0$ we would get a contradiction by comparing the Pohozaev identities of $v_k$ and $V_k$.

Now we use the notation $W_k$ again and use Proposition \ref{key-w8-8} to rewrite the equation for $W_k$.
Let
$$W_k(z)=\hat w_k(e_1+\epsilon_k z), \quad |z|< \delta_0 \epsilon_k^{-1} $$
for $\delta_0>0$ small. Then from Proposition \ref{key-w8-8} we have
\begin{equation}\label{h-exp}
\mathfrak{h}_k(\delta_ky)=\mathfrak{h}_k(\delta_k e_1)+\delta_k \nabla \mathfrak{h}_k(\delta_k e_1)(y-e_1)+O(\delta_k^2)|y-e_1|^2,
\end{equation}
\begin{equation}\label{y-1}
|y|^{2N}=|e_1+\epsilon_k z|^{2N}=1+O(\epsilon_k)|z|,
\end{equation}
\begin{equation}\label{v-radial}
V_k(e_1+\epsilon_k z)+2\log \epsilon_k=U_k(z)+O(\epsilon_k)|z|+O(\epsilon_k^2)(\log (1+|z|))^2
\end{equation}
and
\begin{equation}\label{xi-radial}
\xi_k(e_1+\epsilon_k z)+2\log \epsilon_k=U_k(z)+O(\epsilon_k)(1+|z|).
\end{equation}
Using (\ref{h-exp}),(\ref{y-1}),(\ref{v-radial}) and (\ref{xi-radial}) in (\ref{hat-w}) we write the equation of $W_k$ as
\begin{equation}\label{W-rough}
\Delta W_k+\mathfrak{h}_k(\delta_k e_1)e^{U_k(z)}W_k=-\epsilon_k \nabla \mathfrak{h}_k(0)\cdot ze^{U_k(z)}+E_w, \quad 0<|z|<\delta_0 \epsilon_k^{-1}
\end{equation}
where
\begin{equation}\label{rough-Ew}
E_w=O(\delta_k)(1+|z|)^{-4}+O(\epsilon_k)(1+|z|)^{-3}W_k(z), \quad |z|<\delta_0 \epsilon_k^{-1}.
\end{equation}
At this moment we use $|W_k(z)|\le C$ and a rough estimate of $E_w$ is
\begin{equation}\label{ew-rough-2}
E_w(z)=O(\epsilon_k)(1+|z|)^{-3}, \quad |z|<\delta_0 \epsilon_k^{-1}.
\end{equation}

Since $\hat w_k$ obviously converges to a global harmonic function with removable singularity, we have $\hat w_k\to \bar c$ for some $\bar c\in \mathbb R$. Then we claim that
\begin{lem}\label{bar-c-0} $\bar c=0$.
\end{lem}

\noindent{\bf Proof of Lemma \ref{bar-c-0}:}

 If $\bar c_1\neq 0$, we use $W_k(z)=\bar c+o(1)$ on $B(0,\delta_0 \epsilon_k^{-1})\setminus B(0, \frac 12\delta_0 \epsilon_k^{-1})$ and consider the projection of $W_k$ on $1$:
 $$g_0(r)=\frac 1{2\pi}\int_{0}^{2\pi}W_k(re^{i\theta})d\theta. $$
If we use $F_0$ to denote the projection to $1$ of the right hand side we have, using the rough estimate of $E_w$ in (\ref{ew-rough-2})
$$g_0''(r)+\frac 1r g_0'(r)+\mathfrak{h}_k(\delta_ke_1)e^{U_k(r)}g_0(r)=F_0,\quad 0<r<\delta_0 \epsilon_k^{-1} $$
where
$$F_0(r)=O(\epsilon_k)(1+|z|)^{-3}. $$
In addition we also have
$$\lim_{k\to \infty} g_0(\delta_0 \epsilon_k^{-1})=\bar c_1+o(1). $$
Here the $O(\delta_k)(1+|z|)^{-4}$ is absorbed because of the smallness of $\delta$.
For simplicity we omit $k$ in some notations. By the same argument as in Lemma \ref{w-around-e1},  we have
$$g_0(r)=O(\epsilon_k)(\log (2+r))^2,\quad 0<r<\delta_0 \epsilon_k^{-1}. $$
Thus $\bar c_1=0$.
Lemma \ref{bar-c-0} is established. $\Box$

\medskip

Based on Lemma \ref{bar-c-0} and standard Harnack inequality for elliptic equations we have
\begin{equation}\label{small-til-w}
\tilde w_k(x)=o(1),\,\,\nabla \tilde w_k(x)=o(1),\,\, x\in B_3\setminus (\cup_{l=1}^N(B(e^{i\beta_l},\delta_0)\setminus B(e^{i\beta_l}, \delta_0/8))).
\end{equation}
Equation (\ref{small-til-w}) is equivalent to $w_k=o(\delta_k)$ and $\nabla w_k=o(\delta_k)$ in the same region.

\medskip

\noindent{\bf Proof of Theorem \ref{main-thm}:}

 For $s=1,...,N$ we consider the Pohozaev identity around $Q_s^k$. Let $\Omega_{s,k}=B(Q_s^k,r)$ for small $r>0$. For $v_k$ we have
\begin{align}\label{pi-vk}
\int_{\Omega_{s,k}}\partial_{\xi}(|y|^{2N}\mathfrak{h}_k(\delta_ky))e^{v_k}-\int_{\partial \Omega_{s,k}}e^{v_k}|y|^{2N}\mathfrak{h}_k(\delta_ky)(\xi\cdot \nu)\\
=\int_{\partial \Omega_{s,k}}(\partial_{\nu}v_k\partial_{\xi}v_k-\frac 12|\nabla v_k|^2(\xi\cdot \nu))dS. \nonumber
\end{align}
where $\xi$ is an arbitrary unit vector. Correspondingly the Pohozaev identity for $V_k$ is

\begin{align}\label{pi-Vk}
\int_{\Omega_{s,k}}\partial_{\xi}(|y|^{2N}\mathfrak{h}_k(\delta_ke_1))e^{V_k}-\int_{\partial \Omega_{s,k}}e^{V_k}|y|^{2N}\mathfrak{h}_k(\delta_k e_1)(\xi\cdot \nu)\\
=\int_{\partial \Omega_{s,k}}(\partial_{\nu}V_k\partial_{\xi}V_k-\frac 12|\nabla V_k|^2(\xi\cdot \nu))dS. \nonumber
\end{align}
Using $w_k=v_k-V_k$ and $|w_k(y)|\le C\delta_k$ we have
\begin{align*}
&\int_{\partial \Omega_{s,k}}(\partial_{\nu}v_k\partial_{\xi}v_k-\frac 12|\nabla v_k|^2(\xi\cdot \nu))dS\\
=&\int_{\partial \Omega_{s,k}}(\partial_{\nu}V_k\partial_{\xi}V_k-\frac 12|\nabla V_k|^2(\xi\cdot \nu))dS\\
&+\int_{\partial \Omega_{s,k}}(\partial_{\nu}V_k\partial_{\xi}w_k+\partial_{\nu}w_k\partial_{\xi}V_k-(\nabla V_k\cdot \nabla w_k)(\xi\cdot \nu))dS+O(\delta_k^2).
\end{align*}
If we just use crude estimate: $\nabla w_k=o(\delta_k)$, then
$$\int_{\partial \Omega_{s,k}}(\partial_{\nu}v_k\partial_{\xi}v_k-\frac 12|\nabla v_k|^2(\xi\cdot \nu))dS-
\int_{\partial \Omega_{s,k}}(\partial_{\nu}V_k\partial_{\xi}V_k-\frac 12|\nabla V_k|^2(\xi\cdot \nu))dS=o(\delta_k).
$$
The difference on the second terms is minor:
$$\int_{\partial \Omega_{s,k}}e^{v_k}|y|^{2N}\mathfrak{h}_k(\delta_ky)(\xi\cdot \nu)-\int_{\partial \Omega_{s,k}}e^{V_k}|y|^{2N}\mathfrak{h}_k(\delta_ke_1)(\xi\cdot \nu)=O(\delta_k \epsilon_k^2). $$
To evaluate the first term,  we use
\begin{align}\label{imp-1}
&\partial_{\xi}(|y|^{2N}\mathfrak{h}_k(\delta_ky))e^{v_k}\\
=&\partial_{\xi}(|y|^{2N}\mathfrak{h}_k(\delta_ke_1)+|y|^{2N}\delta_k\nabla \mathfrak{h}_k(\delta_ke_1)(y-e_1)+O(\delta_k^2))e^{V_k}(1+w_k+O(\delta_k^2))\nonumber\\
=&\partial_{\xi}(|y|^{2N})\mathfrak{h}_k(\delta_k e_1)e^{V_k}+\delta_k\partial_{\xi}(|y|^{2N}\nabla \mathfrak{h}_k(\delta_ke_1)(y-e_1))e^{V_k}\nonumber\\
&+\partial_{\xi}(|y|^{2N}\mathfrak{h}_k(\delta_ke_1))e^{V_k}w_k+O(\delta_k^2)e^{V_k}.\nonumber
\end{align}

For the third term on the right hand side of (\ref{imp-1}) we use the equation for $w_k$:
$$\Delta w_k+\mathfrak{h}_k(\delta_ke_1)e^{V_k}|y|^{2N}w_k=-\delta_k\nabla \mathfrak{h}_k(\delta_ke_1)\cdot (y-e_1)|y|^{2N}e^{V_k}+O(\delta_k^2)e^{V_k}|y|^{2N}.$$
From integration by parts we have
\begin{align}\label{extra-1}
&\int_{\Omega_{s,k}}\partial_{\xi}(|y|^{2N})\mathfrak{h}_k(\delta_ke_1)e^{V_k}w_k\nonumber\\
=&2N\int_{\Omega_{s,k}}|y|^{2N-2}y_{\xi}\mathfrak{h}_k(\delta_ke_1)e^{V_k}w_k\nonumber\\
=&2N\int_{\Omega_{s,k}}\frac{y_{\xi}}{|y|^2}(-\Delta w_k-\delta_k \nabla\mathfrak{h}_k(\delta_ke_1)(y-e_1)|y|^{2N}e^{V_k}+O(\delta_k^2)e^{V_k})\nonumber\\
=&-2N\delta_k\int_{\Omega_{s,k}}\frac{y_{\xi}}{|y|^{2}}\nabla \mathfrak{h}_k(\delta_ke_1)(y-e_1)|y|^{2N}e^{V_k}\nonumber\\
&+2N\int_{\partial \Omega_{s,k}}(\partial_{\nu}(\frac{y_{\xi}}{|y|^2})w_k-\partial_{\nu}w_k\frac{y_{\xi}}{|y|^2})+O(\delta_k^2)\nonumber\\
=&-16N\delta_k\pi(e^{i\beta_s}\cdot \xi)\nabla\mathfrak{h}_k(\delta_ke_1)(e^{i\beta_s}-e_1)+o(\delta_k),
\end{align}
where we have used $\nabla w_k, w_k=o(\delta_k)$ on $\partial \Omega_{s,k}$.
For the second term on the right hand side of (\ref{imp-1}), we have
\begin{align}\label{imp-2}
&\int_{\Omega_{s,k}}\delta_k\partial_{\xi}(|y|^{2N}\nabla \mathfrak{h}_k(\delta_ke_1)(y-e_1))e^{V_k}\\
=&2N\delta_k\int_{\Omega_{s,k}}y_{\xi}|y|^{2N-2}\nabla \mathfrak{h}_k(\delta_ke_1)(y-e_1)e^{V_k}+
\delta_k\int_{\Omega_{s,k}}|y|^{2N}\partial_{\xi}\mathfrak{h}_k(\delta_ke_1)e^{V_k} \nonumber \\
=&16N\pi\delta_k(e^{i\beta_s}\cdot \xi)\nabla \mathfrak{h}_k(\delta_ke_1)(e^{i\beta_s}-e_1)+
8\pi\delta_k\partial_{\xi}\mathfrak{h}_k(\delta_ke_1)+o(\delta_k).\nonumber
\end{align}
Using (\ref{extra-1}) and (\ref{imp-2}) in the difference between (\ref{pi-vk}) and (\ref{pi-Vk}), we have
$$\delta_k \partial_{\xi}\mathfrak{h}_k(\delta_ke_1)=o(\delta_k). $$
Thus $\nabla \mathfrak{h}_k(\delta_ke_1)=o(1)$. Theorem \ref{main-thm} is established. $\Box$

\end{document}